\documentclass[preprint,12pt]{elsarticle}

\usepackage{amssymb,amsmath,bm}
\usepackage{graphicx}
\newcommand{\xx}{\mathbf{x}}
\newcommand{\CC}{\mathbf{C}}
\newcommand{\Cov}{\mathbf{\Gamma}}

\begin{document}

\begin{frontmatter}



\title{Innovative observing strategy and orbit determination for Low
  Earth Orbit Space Debris}


\author[dm]{A. Milani}
\author[dm]{D. Farnocchia}
\author[dm]{L. Dimare}
\author[cnr]{A. Rossi}
\author[dm]{F. Bernardi}

\address[dm]{Department of Mathematics, University of Pisa, Largo
  Bruno Pontecorvo 5, 56127, Pisa, Italy} 

\address[cnr]{IFAC-CNR \& ISTI-CNR, CNR Area della Ricerca di Firenze,
  Via Madonna del Piano 10, 50019, Sesto Fiorentino, Firenze, Italy}

\begin{abstract}
  We present the results of a large scale simulation, reproducing the
  behavior of a data center for the build-up and maintenance of a
  complete catalog of space debris in the upper part of the low Earth
  orbits region (LEO). The purpose is to determine the performances of
  a network of advanced optical sensors, through the use of the newest
  orbit determination algorithms developed by the Department of
  Mathematics of Pisa (DM). Such a network has been proposed to ESA in
  the Space Situational Awareness (SSA) framework by Carlo Gavazzi
  Space SpA (CGS), Istituto Nazionale di Astrofisica (INAF), DM, and
  Istituto di Scienza e Tecnologie dell'Informazione (ISTI-CNR).  The
  conclusion is that it is possible to use a network of optical
  sensors to build up a catalog containing more than 98\% of the
  objects with perigee height between 1100 and 2000 km, which would be
  observable by a reference radar system selected as comparison. It is
  also possible to maintain such a catalog within the accuracy
  requirements motivated by collision avoidance, and to detect
  catastrophic fragmentation events. However, such results depend upon
  specific assumptions on the sensor and on the software technologies.
\end{abstract}


\end{frontmatter}


\section{Introduction}
In the context of the European program SSA, the aim of the project
\textit{SARA-Part I Feasibility study of an innovative system for
  debris surveillance in LEO regime} was to demonstrate the
feasibility of a European network based on optical sensors, capable of
complementing the use of radars for the identification and
cataloging of debris in the high part of the LEO region, to lower the
requirements on the radar system in terms of power and
performances. The proposal relied on the definition of a wide-eye
optical instrument able to work in high LEO zone and on the
development of new orbit determination algorithms, suitable for the
kind and amount of data coming from the surveys targeting at LEOs.

Taking into account the performances expected from the innovative optical
sensor, we have been able to define an observing strategy allowing to
acquire data from every object passing above a station of the assumed
network, provided the $S/N$ is good enough. Still the number of
telescopes required for a survey capable of a rapid debris catalog
build-up was large: to reduce this number we have assumed that the goal
of the survey had to be only one exposure per pass.

A new algorithm, based on the first integrals of the Kepler problem,
was developed by DM to solve the critical issue of LEOs orbit
determination. Standard methods, such as Gauss' \cite{gauss}, require
at least three observations per pass in order to compute a preliminary
orbit, while the proposed algorithm needs only two exposures, observed
at different passes. This results in a significant reduction of the
number of required telescopes, thus of the cost of the entire
system. For LEO, the proposed method takes into account the nodal
precession due to the quadrupole term of the Earth
geopotential. Because of the low altitude of the orbits and the
availability of sparse observations, separated by several orbital
periods, this effect is not negligible and it must be considered since
the first step of preliminary orbit computation.

The aim was to perform a realistic simulation. Thus in addition to the
correlation and orbit determination algorithms, all the relevant
elements of the optic system were considered: the telescope design,
the network of sensors, the observation constraints and strategy, the
image processing techniques. Starting from the ESA-MASTER2005
population model, CGS provided us with simulated observations,
produced taking into account the performances of the optical
sensors. These data were processed with our new orbit determination
algorithms in the simulations of three different operational phases:
catalog build-up, orbit improvement, and fragmentation analysis. The
results of these simulations are given in the
Sections~\ref{sec:survey_res}, \ref{sec:task_results}, and
\ref{sec:frag}. 

\section{Assumptions}
\label{sec:assumptions}

The only way to validate a proposed system, including a network of
sensors and the data processing algorithms, was to perform a realistic
simulation. This does not mean a simulation including all details, but
one in which the main difficulties of the task are addressed.  The
output of such a simulation depends upon all the assumptions used, be
they on the sensor performance, on the algorithms and software, on the
physical constraints (e.g., meteorological conditions).

In what follows, we list all the assumptions used in the catalog
build-up simulation, in the orbit improvement simulation, and in the
fragmentation detection simulation.  We discuss the importance of each
one in either surveying capability or detection capability or orbit
availability and accuracy. All the assumptions turn out to be
essential to achieve the performance measured by the simulations.

\subsection{Assumptions on Sensors}

We are assuming a network consisting of optical sensors only, with the
following properties:

\begin{enumerate}

\item The telescope and the camera have a \textbf{large field of view:
    45 square degrees} ($6.66^\circ\times 6.66^\circ =24000\times 24000$ arcsec).
\item The telescope has \textbf{quick motion capability}, with
  mechanical components allowing a \textbf{1 s exposure every 3 s},
  with each image covering new sky area: the motion in the 2 s
  interval must be $\ge 6.66^\circ$, with stabilization in the new position.
\item The camera system has a \textbf{quick readout}, to be able to
  acquire the image from all the CCD chips \textbf{within the same 2
    s} used for telescope repositioning, and this with a low readout
  noise, such that the main source of noise is the sky background.
\item The camera system needs to have \textbf{high resolution},
  comparable to the seeing. A pixel scale of about 1.5 arcsec is the
  best compromise, known to allow for accurate astrometry. Then field
  of view of 24000$\times$24000 arcsec implies a camera system with
  256 MegaPixel.
\item The camera has by design a \textbf{fill factor 1}.  The fill
  factor is the ratio between the effective area, on the focal plane,
  of the active sending elements and the area of the field of view.
\item The \textbf{telescope aperture needs to be large} enough to
  detect the target debris population, we are assuming an
  \textbf{effective aperture of 1 meter}. That is, the unobstructed
  photon collecting area has an area equal to a disk of 1 meter
  diameter.
\item The \textbf{network of sensors} includes \textbf{7
    geographically distributed stations, each with 3 telescopes}
  available for LEO tracking.
\item The telescope is assumed to have \textbf{tasking capabilities},
  consisting in the possibility of \textbf{non-sidereal tracking at a
    programmed rate up to 2000 arcsec/s} (relative to the sidereal
  frame) while maintaining the image stable.
\end{enumerate}

The above assumptions about the sensor hardware require a significant
effort in both technological development and resources: a discussion
on the feasibility of each of them is beside the scope of this
paper. We need just to point out that a design for an innovative
sensor with such properties does exist and has been presented in
\cite{cibin}.

The large field of view is needed to cope with the tight requirements
on surveying capability resulting from a population of space objects
with very fast angular motion, up to 2000 arcsec/s. The surveying
capability is further enhanced by taking an image on a new field of
view every 3 s, and by the use of 21 telescopes. Thus it should not be
surprising that such a network has the capability of observing objects
in orbits lower than those previously considered suitable for optical
tracking.

The large fill factor also contributes by making the surveying
capability deterministic, that is objects in the field of view are
effectively observable every time. Manufacturing imperfections
unavoidably decrease the fill factor, but with good quality chips the
reduction to a value around $0.98$ does not substantially change the
performance, while detectors with values in the $0.7\div 0.8$ range
would result in a severely decreased performance.

Both properties (large aperture and large fill factor) are feasible as
a result of an innovative design based on the fly-eye concept, that is
the telescope does not have a monolithic focal plane but many of them,
each filled by the sensitive area of a separate, single chip camera.

The comparatively high astrometric resolution is essential to
guarantee the observation accuracy, in turn guaranteeing the orbit
determination accuracy. If the goal is to perform collision avoidance,
there is no point in having low accuracy orbits. Thus a low
astrometric resolution survey would require a separate tasking/follow
up network of telescopes for orbit improvement. In our assumed network
all the tasking is performed with the same telescopes.

The telescope effective aperture, which corresponds in the proposed
design to a primary mirror with $1.1$ meters diameter, is enough to
observe LEOs in the $8\div 10$ cm diameter range if it is coupled with
a computationally aggressive image processing algorithm, discussed
below and in Sec.~\ref{sec:s2n}.

The selection of the stations locations is a complicated problem,
because it has to strike a compromise between the requirement of a
wide geographical distribution and the constraints from meteorology,
logistics and geopolitics; see Sec.~\ref{sec:network}. To simulate the
outcome of such a complicated selection process, we have used a
network which is ideal neither from the point of view of geographical
distribution nor for meteorological conditions, but it is quite
realistic.

The assumed sensors are optimized for LEO, but they are also very
efficient to observe Medium Earth Orbit (MEO), Geostationary Orbit
(GEO), and any other Earth orbit above 1000 km. The same sensors could
also be used for Near Earth Objects; the required changes have only to
do with longer exposure times and can be implemented in software.

\subsection{Data processing assumptions}

We are assuming the observations from the optical sensor network are
processed with algorithms and the corresponding software, having the
following properties:

\begin{enumerate}
\setcounter{enumi}{8}
\item The \textbf{scheduler} of the optical observations is capable of
  \textbf{taking into account the geometry of light and the phase}
  (which is defined as the Sun-object-observer angle), in such a way
  that the objects passing above the station are imaged and the phase
  is minimized.
\item The \textbf{image processing} includes a procedure to
  \textbf{detect long trails at low signal to noise}, $S/N$ with a
  loss due to the spreading of the image on $T$ pixels proportional to
  $\sqrt{T}$.
\item The \textbf{astrometric reduction algorithms} allow for
  \textbf{sub-pixel accuracy, even for long trails}, and taking
  properly into account star catalog errors.
\item The \textbf{correlation and orbit determination} algorithms
  allow us to compute preliminary orbits starting from a \textbf{single
    trail per pass}, and correlating passes separated by several
  orbital periods of the objects (e.g., a time span of the order of a
  day).
\end{enumerate}

The assumptions 9-12 are different from the previous ones because they
are all about software. Of course a significant software development
effort is necessary, but we consider part of our current research
effort to ensure that the algorithms which could lead to the assumed
results already exist.


The synthetic observations used in the simulation have been obtained
by taking one exposure for each pass, in the visibility interval ($>$
15$^\circ$ of elevation, illuminated by the Sun, station in
darkness). Within this interval, we have assumed the best third from
the point of view of the phase angle is used (near the shadow of the
Earth). Since the apparent magnitude $h$ of the objects is a steep
function of the phase angle, the number of observations with
sufficient $S/N$ is significantly increased (by a factor 3-4).  Such a
\textit{light aware} scheduler was not actually available, but we have
tested that a simple observing strategy exists leading to this
result. The idea is to use a dynamic barrier formed by frequently
visited fields of view. The barrier could be bordering the Earth
shadow at the altitude of the objects being targeted: both simple
computations and a numerical simulation show that this can be achieved
by using $2\div 3$ telescopes with the performances outlined in
assumptions 1-6.


The \textit{trailing loss}, that is the decrease of the signal due to
the spreading of the image on $T$ pixels, appears to limit the
sensitivity of the detector for objects with a high angular velocity,
like $300 \div 1000$ arcsec/s (typical values for an object at an
altitude of 1400 km). This appears to defeat the approach used in
astronomy, that is increasing the exposure time to observe dimmer
objects. However, even for a stationary object such as a star, the
increase in $S/N$ is only with the square root of the exposure
time. Thus, if an algorithms is available to \textit{sum up} the
signal from adjacent pixels, in such a way that $S/N$ accumulates with
$\sqrt{T}$, the increase of exposure time is as effective as for a
stationary target. Such algorithms exist and are discussed in
Sec.~\ref{sec:s2n}. The actual implementation in operational software
and field testing are assumptions.


The observations have to be reduced astrometrically in an accurate
way, with RMS error of 0.4 arcsec when the pixel $S/N$ is good. GEO
and Geostationary Transfer Orbit (GTO) data from the ESA Optical
ground Sensor in Teide (Canary Islands), reduced by University of
Bern, show a typical RMS $0.5\div 0.6$ arcsec of the residuals from
the orbit determination performed by our group \cite{GEO}. For low
$S/N$ on each pixel the astrometric error is assumed to increase, see
Sec.~\ref{sec:s2n}.  When the RMS grows to $>2$ arcsec, the
observations can result in orbit determination failure and/or accuracy
requirements non compliance.
Improvements in the astrometric reduction procedure, to remove
systematic star catalog errors, are already implemented in asteroid
orbit catalogs such as the ones available online from the AstDyS-2
and NEODyS-2 development systems. They are based on the star catalog
debiasing algorithms proposed by \cite{cbm10}.  The assumption is that
an ad hoc astrometric reduction software is developed.


Correlation and Orbit Determination algorithms, developed and
implemented in software by our group, have the capability to use
significantly less observations with respect to classical methods to
compute a preliminary orbit. As an example, we can use two trails from
different passes of the same object above either the same or a
different station to compute an orbit with covariance matrix. These
methods, and those used for successive orbit improvement, are
discussed in Sec.~\ref{sec:orbdet}.
The correlation software we use is capable of computing orbits
starting from sparse uncorrelated observations of LEO (also of GEO,
\cite{GEO}). The amount of data to be used as input for initial
catalog build-up is limited to 1 exposure per pass.  The advantage with
respect to the traditional approach \cite{gauss,escobal}, requiring
three separate observations in the same pass, is such that the
surveying capability for a given sensor network is increased by a
factor 3.
Note that in this case the software with the assumed performances
actually exists and is being tested in simulations like the ones we
are discussing in this paper. The assumption is only that the existing
software is upgraded to operational.

\section{Observing Strategy}
\label{sec:obs_strategy}

To perform a debris observation some conditions shall be verified: a
minimum elevation angle, the orbiting object must be in sunlight,
etc.. These conditions are strongly dependent on the object orbit
parameters, on the observatory location and on the seasonal
factors. There are also other observational constraints that have been
taken into account, such as the distance from the Moon and the
galactic plane.

\subsection{Geometrical constraints}
\label{sec:horizon}

\begin{figure}[h!]
\centering
\includegraphics[width=0.7\textwidth]{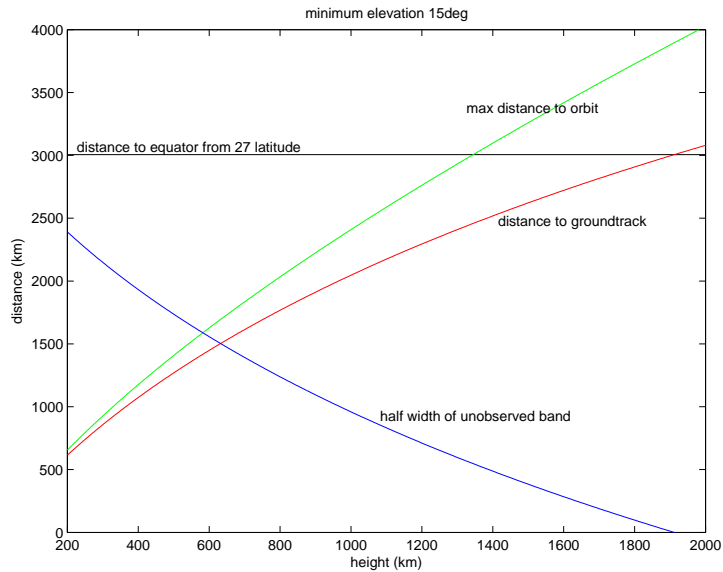}
\caption{Maximum distance to orbit, to groundtrack and half
  width of the equatorial unobservable belt}
\label{fig:horizon}
\end{figure}

The first constraints to the network architecture are purely
geometrical and are due to the horizon. An orbiting object at an
altitude $h_p$ is visible only up to a given distance from a station,
beyond which the object is below a minimum elevation, $15^\circ$ being
a reasonable value. For an object at $h=1400$ km, the distance to the
object is thus limited to about 3100 km, and the distance to the
groundtrack of the object to about 2500 km; see
Fig.~\ref{fig:horizon}, showing these values in km as a function of
the object altitude. Moreover, for a station at a latitude of
$27^\circ$, this Figure also shows the half width of the equatorial
band such that, if the groundtrack is in there, the object is not
observable. This argument favors the stations located at low
latitudes.

\begin{figure}[h!]
\centering
\includegraphics[width=0.7\textwidth]{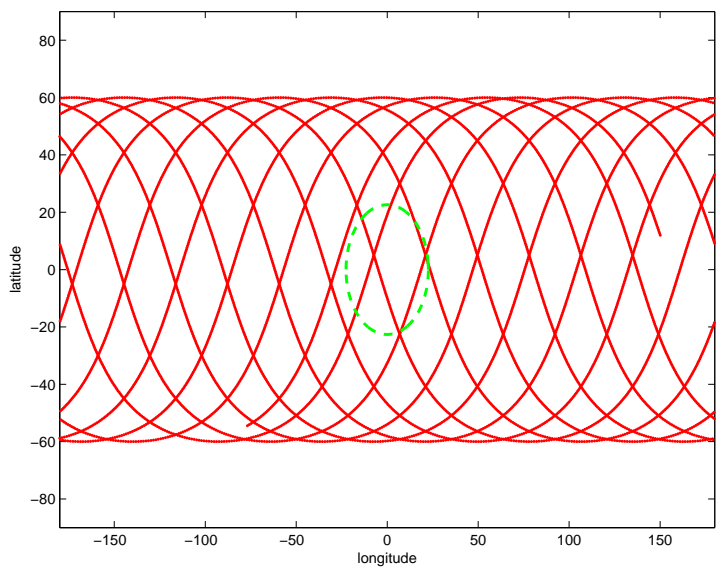}
\caption{Ground track for $\mathbf{h_p}=$ 1400 km,
  $I=$60$^{\mathbf \circ}$}
\label{fig:groundtra}
\end{figure}

The second consideration is how the object presents itself with a
groundtrack passing near a station.  Figure~\ref{fig:groundtra} shows
the groundtrack for a nearly circular orbit with $h_p=1400$ km and
inclination of 60$^\circ$. The oval contour shows the maximum
visibility range, for an equatorial station and an elevation $\ge
15^\circ$.  Typically a LEO has 4 passes/day above the required
altitude as seen from a station at low latitudes. Note that the
constraints discussed so far apply equally to a radar sensor.

\subsection{Geometry of sunlight}
\label{sec:geomlight}
 
The main difference with radar arises
because of the geometry of sunlight. The requirements for an optical station
are the following:
\begin{enumerate}
\item the ground station is in darkness, e.g., the Sun must be at
  least 10-12$^\circ$ below the horizon, that is the sky is dark
  enough to begin operations, typically about 30-60 minutes after
  sunset and before sunrise (this is strongly dependent on the latitude
  and the season at the station);
\item the orbiting object is in sunlight; 
\item the atmosphere is clear (no dense clouds).
\end{enumerate}

\begin{figure}[t]
\centering
\includegraphics[width=0.8\textwidth]{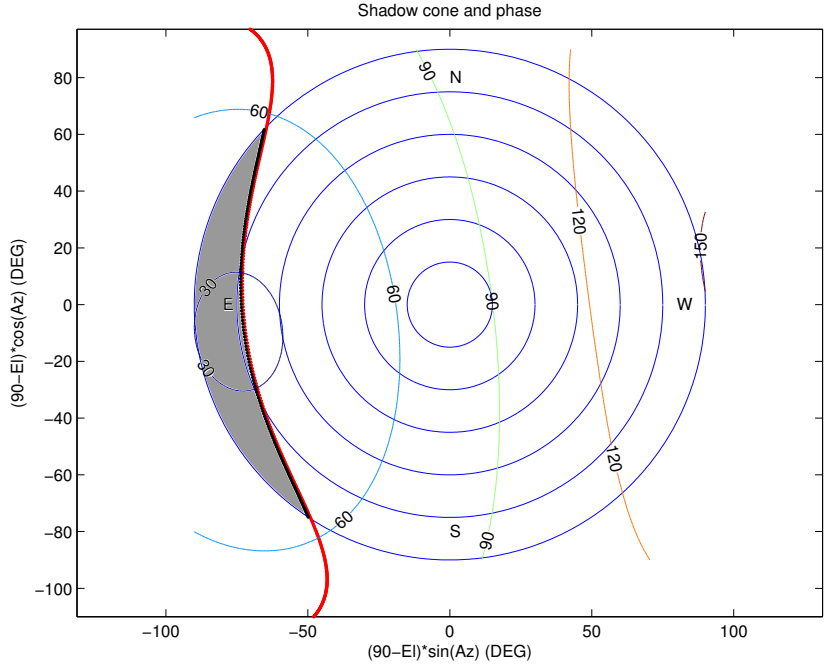}
\caption{Earth Shadow and iso-phase curves, tropical station at solar
  time 19 hours}
\label{fig:shadow1}
\end{figure}

\begin{figure}[t]
\centering
\includegraphics[width=0.8\textwidth]{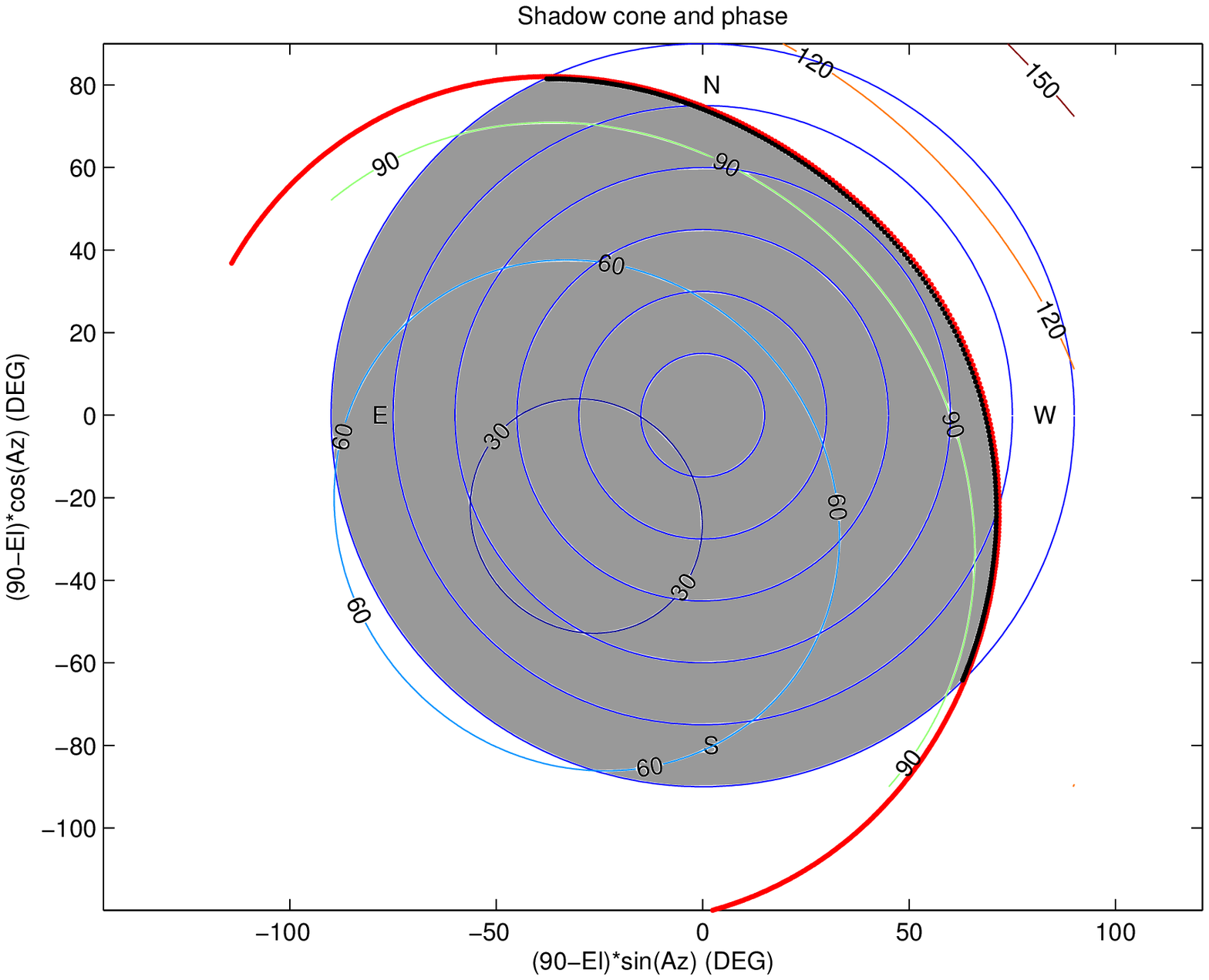}
\caption{Earth shadow and iso-phase curves, tropical station at solar
  time 22 hours}
\label{fig:shadow2}
\end{figure}

\begin{figure}[t]
\centering
\includegraphics[width=0.8\textwidth]{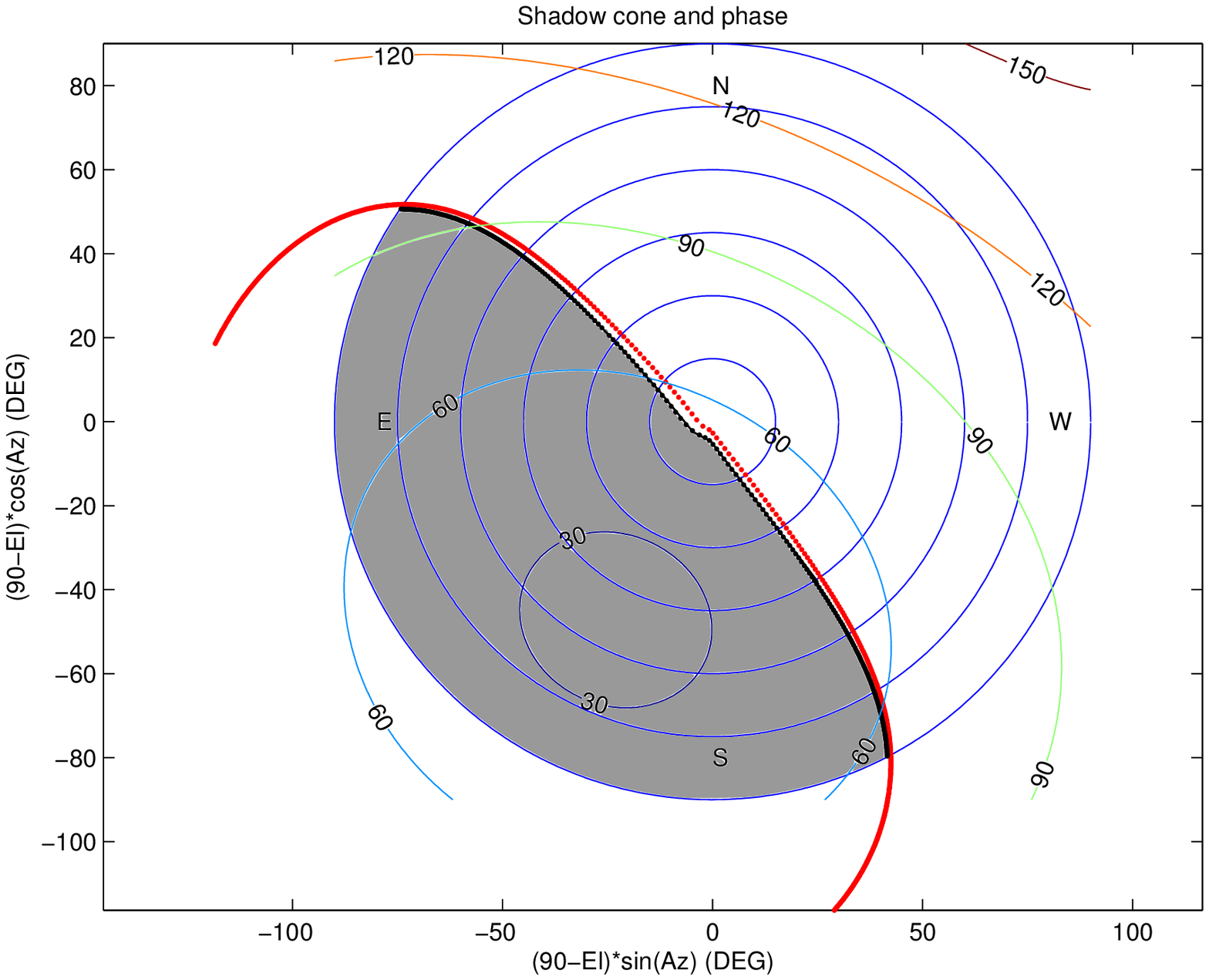}
\caption{Earth shadow and iso-phase curves, high latitude
  station at solar time 22 hours}
\label{fig:shadow3}
\end{figure}

The condition on sunlight is quite restrictive: the low orbiting
objects are fully illuminated in all directions only just after
sunset and before sunrise. In the
Figures~\ref{fig:shadow1}-\ref{fig:shadow3} the bold line represents the
Earth shadow boundary at 1400 km above ground. The shadow also depends
upon the season, the figures have been drawn for March 20, a date
close to an equinox. The Earth shadow region, where the orbiting
object is invisible, is represented in gray. The circles represent the
iso-elevation regions of the sky above the horizon, which is the outer
curve; the center is the local zenith.


The labeled lines (30, 60, 90 and 120) represent the iso-phase
curves for objects at 1400 km above ground, that is the directions in
the sky where the objects have a specific phase angle.
The phase angle is a very critical observing parameter for a
debris. The optical
magnitude of an object (generally all Solar System moving objects)
depends, among other parameters, by the phase angle: the smaller the
phase angle the brighter the object. The strength of this effect also
depends upon the optical properties of the object surface, such as the
albedo. Anyway the effect is large, e.g., at a phase of $90^\circ$ the
apparent magnitude could increase by $>3$ magnitudes with respect to
an object at the same distance but with $0^\circ$ phase. Thus an
observation scheduling taking into account the need of observing with
the lowest possible phase increases very significantly the optical
sensor performance.

The Figures show that the regions where the phase angles are smaller
are close to the Earth shadow boundary. Very low phases can be
achieved only near sunset and sunrise, by looking in a direction
roughly opposite to the Sun (see Fig.~\ref{fig:shadow1}).  For LEOs
there is a central portion of the night, lasting several hours, in
which for a either an equatorial or a tropical station the
observations are either impossible or with a very unfavorable phase,
e.g., for the station at $27^\circ$ North latitude on March 20 from about
22 hours to 2 hours of the next day (see Fig.~\ref{fig:shadow2}). On
the contrary, for a station at an high latitude (both North and South) the
dark period around midnight when LEO cannot be observed does not
occur, because the Earth shadow moves South (for a North station; see
Fig.~\ref{fig:shadow3}).

\subsection{Selection of the observing network}
\label{sec:network}

It is clear that, by combining the orbital geometry of passages above
the station with the no shadow condition, it is possible to obtain
objects which are unobservable from any given low latitude station, at
least for a time span until the precession of the orbit (due to
Earth's oblateness, $\sim 5^\circ\, \cos I$/day for an altitude of
1400 km and inclination $I$) changes the angle between the orbit plane
and the direction to the Sun. On the other hand, high latitude
stations cannot observe low inclination objects, and operate for a
lower number of hours per year because of shorter hours of darkness in
summer and worse weather in winter. A trade off is needed, which
suggests to select some intermediate latitude stations, somewhere
between 40 and 50$^\circ$ both North and South.

The meteorological constraints can be handled by having multiple
opportunities of observations from stations far enough to have low
meteorological correlation. This implies that an optimal network needs
to include both tropical and high latitude stations, with a good
distribution also in longitude.
Beside the need for geographic distribution as discussed above, the
other elements to be considered in the selection of the network are
the following:
\begin{itemize}
\item Geopolitics: the land needs to belong to Europe, or to friendly
  nations. The limitations
  of the European continent implies that sites in minor European
  islands around the world are needed as well as observing sites in
  other countries.
\item Logistics: some essentials like electrical power, water supply,
  telecommunications, airports, harbors, and roads have to be available.
\item Meteo: the cloud cover can be extremely high, in some geographic
  areas, and especially in (local) winter. Other meteorological
  parameters such as humidity, seeing, wind play an important
  role. High elevation observing sites over the inversion layer are
  desirable, but in mid ocean there are not many mountains high enough.
\item Orography: an unobstructed view of the needed sky portions, down
  to 15$^\circ$ of elevation, is necessary. Astronomical observatories
  are not so demanding, especially in the pole direction.
\item Light pollution: an observing site with low light pollution is
  essential, to lower the sky background, which is the main source of
  noise.
\end{itemize}

\begin{table}
\small
\caption{\label{tab:network} Geographical coordinates of the proposed 
  network of stations}
\centering
\begin{tabular}{|r|r|r|r|}
\hline 
STATION & Latitude [deg]& Longitude [deg]& Height [m]\\
\hline
TEIDE  & 28$^\circ$ 18' 03.3'' N & 16$^\circ$ 30' 42.5'' W & 2390\\
(Canary Islands, Spain)&&&\\
\hline
HAO  & 18$^\circ$ 08' 45.0'' S & 140$^\circ$ 52' 54.0'' W &  6\\
(French Polynesia, France)&&&\\
\hline
FALKLAND ISLAND & 51$^\circ$ 42' 32.0'' S & 57$^\circ$ 50' 27.0'' W & 30\\
(Great Britain) &&&\\
\hline
NEW NORCIA  & 31$^\circ$ 02' 54.0'' S & 116$^\circ$ 11' 31.0'' E & 244\\
(Australia) &&&\\
\hline
MALARGUE  & 35$^\circ$ 46' 24.0'' S & 69$^\circ$ 23' 59.0'' W & 1509\\
(Argentina) &&&\\
\hline
GRAN SASSO  & 42$^\circ$ 29' 60.0 N & 13$^\circ$ 33' 04.0'' E & 1439\\
(Italy) &&&\\
\hline
PICO DE VARA  & 37$^\circ$ 47' 48.0'' N & 25$^\circ$ 13' 10.0'' W & 579\\
(Azores Islands, Portugal) &&&\\
\hline
\end{tabular}
\end{table}

The optical station network used for our simulations takes into
account, as much as possible, the above constraints and is shown on
Table
. In the selection of the stations we took into
account the local meteorological statistics. This was done by
accessing meteorological data bases, mostly data obtained from the
ISCCP (NASA) project (http://isccp.giss.nasa.gov/index.html). The
analysis of these data forced us to give up the possibility of using
some geographically convenient locations, which turned out to have a
``total cloud amount'', defined as the percentage of cloud coverage of
the sky for each day of the year, far too high: e.g., the islands of
St. Pierre and Michelon, French territory near Canada, were excluded.

\section{Population model}

In order to produce a realistic simulation of the whole observation
process a suitable population of orbiting objects is required.  A
subset of the ESA MASTER-2005 population model \cite{master}, upgraded
with the recent in-orbit collisions (FengYun-1C, IRIDIUM 33 - COSMOS
2251) was provided by ESOC.  The MASTER model contains the largest
objects taken from the USSTRATCOM Two Line Elements (TLE) plus smaller
objects generated with ad hoc source models. The subset available for
this work included 31686 entries, representing either single objects
or sampled ones, with diameter $d$, 3 cm $< d<$ 31.7 m and crossing
the Low Earth Orbit (LEO) region (i.e., with the perigee altitude
$h_p$ between 200 km and 2000 km).
%

The population file did not contain any value for the albedo of the
objects, a quantity that is needed to derive the magnitude of the
object in the sky.  A commonly accepted value of the albedo for a
generic spacecraft is between 0.1 and 0.2 \cite{africano,kessler}. We
used the conservative assumption of albedo 0.1 for all the objects
considered.  Then the absolute magnitude $H$ was derived according to the
IAU standard for asteroids: $ H=33-5\log_{10}(d)$ where $d$ (in m)
is provided by the MASTER.

\section{Conditions for detection}
\label{sec:s2n}

The possibility to detect an object with an optical sensor depends
upon the $S/N$ of the corresponding image. The $S/N$ is a function of the
following parameters of the object: $H$, distance, and
angular velocity with respect to the image reference frame. The images
are taken in a sidereal reference frame, defined by stars. (Only for
GEOs there is advantage in taking the images in a
reference frame body-fixed to the Earth). The apparent magnitude $h$ can
be computed as follows, if the exposure time is such that there is no
trailing effect (that is, the image has the same shape as the one of a
fixed star):
\begin{equation}
h=H-5 \log_{10} (r) + f(\phi,G)
\end{equation}
where $r$ is the distance station-object in Astronomical Units, $f$ is
the function describing the phase effect, containing the phase angle
$\phi$ and a parameter $G$ depending upon the optical properties of
the surface. As mentioned above, the correction $f$ is important.
The other parameters depend upon the instrument and the atmosphere.
The standard equation to compute $S/N$ for a stationary source
(a star) is given by the following formula:
\[
\left(\frac{S}{N}\right)_{star} = \frac{S}{\sqrt{S+N}} \, .
\]
In this formula $N$ includes the contributions from all the
sources of noise occurring in the measurement, in particular:
\begin{itemize}
\item Read/Out noise: $\sigma^2_{r/o} = n_{pix} \,RN^2 $ (RN =
  pixel r/o noise);
\item Dark current noise: $\sigma^2_{Dk} = D \, n_{pix}\, t$
  (D = dark current);
\item Sky background: $\sigma^2_{sky} = R_{sky}\, t$ (R$_{sky}$ =
  sky background flux).
\end{itemize}
The term $S$ under square root accounts for the Poisson statistics:
$\sigma^2_{Poisson}=S$.

The main problem with optical observations of Earth orbiting objects,
especially for LEO, is their large angular velocity in the sidereal
frame. Thus, unless the exposure times are very short, severely
limiting the maximum magnitude of a detectable object, the image is a
\textit{trail}, spread over a comparatively large number of pixels.
The $S/N$ on a single pixel is given by dividing the total signal of the
debris by $T$; this effect is called \textit{trailing loss}. In other
words, the signal per pixel is $S/T$, thus on a single pixel:
\[
\left(\frac{S}{N}\right)_{pixel} = \frac{{S}/{T}}{\sqrt{{S}/{T}+N}}\ .
\]
The values of $S$ as a function of $h$, and of $N$ as a
function of the instrument properties and of the atmospheric
properties, have been computed on the basis of the telescope design as
in \cite{cibin}. 

\begin{figure}[h!]
\centering
\includegraphics[width=0.8\textwidth]{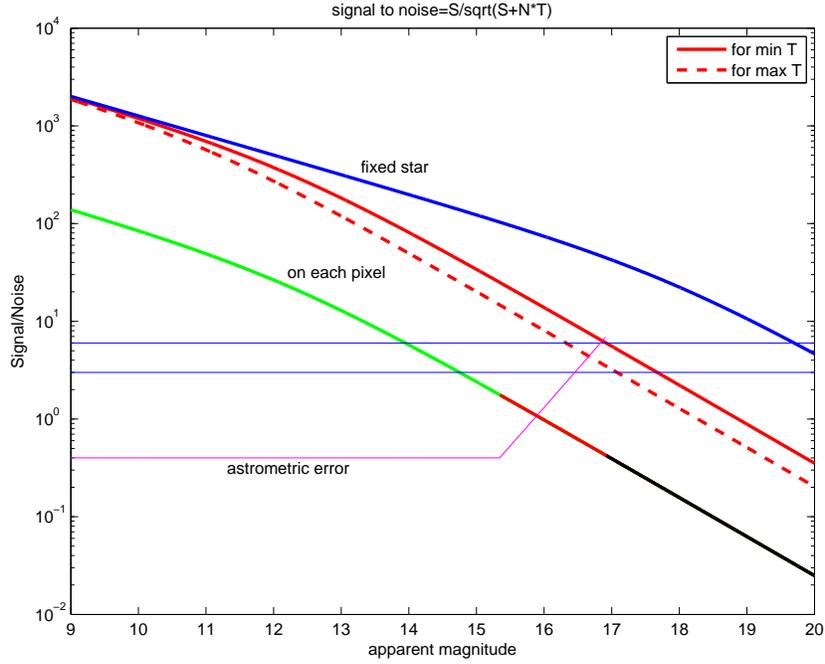}
\caption{$h$ vs. $S/N$, for different trailing loss: from
  the top, $T$=1, 200, 600}
\label{fig:s2n}
\end{figure}

An algorithm to combine the information contained in all the pixels
touched by the trail has been proposed by \cite{villani}. The
principle is to test all the possible trails which could occur in the
frame, with a computationally efficient algorithm to decrease the
computational complexity. Such an algorithm is capable to detect very
faint trails, because the signal along the trail direction is added,
while the noise is added in the RMS sense. Thus the $S/N$ for the trail
is :
\begin{equation}
\label{eq:snr}
\left(\frac{S}{N}\right)_{trail} = \frac{(S/T)
  T}{\sqrt{T}\sqrt{S/T+N}}=\frac{S}{\sqrt{S+NT}}\ .
\end{equation}
Then the $S/N$ of the trail is much larger than the single pixel one
by a factor close to $\sqrt{T}$ for low $S$. According to the tests
reported in \cite{villani}, it is necessary to keep some margin to
avoid false detections of trails, thus we use as criterion
${S}/{N}_{trail} \ge 6$, while for a stationary object lower values,
such as 3, could be used. E.g., for a trail $T=200$ pixels long, the
advantage is by a factor $\sqrt{200}/2\sim 7$, equivalent to more than
2 magnitudes. Fig.~\ref{fig:s2n} shows the $S/N$, as a function of the
$h$, for different observing conditions, resulting in
different values of the trailing loss parameter $T$. From the top: for
a fixed star; for angular velocity $300$ arcsec/s; for $1000$
arcsec/s; on each individual pixel (for $300$ arcsec/s).

The problem is that it is not easy to determine the beginning and the
end of the trail with high accuracy when the trail is too faint. For
faint trails, the astrometric error is determined by the error in
finding the ends, thus is given by:
\begin{equation}
  Z \left(\frac{S}{N}\right)_{pixel} \le \sqrt{Z} \Longrightarrow Z=\frac{TS 
    + T^2N}{S^2}=\left( \frac{1}{(S/N)_{pixel}}\right)^2\ .
\label{eq:astromerr}
\end{equation}

In practice there are two regimes for the astrometric error:
\begin{itemize}
\item when the signal is strong, the astrometric error is dominated by
  the astrometry method, that is by the systematics in the astrometric
  catalogs;
\item when the signal is weak, the astrometric error is determined by
  (\ref{eq:astromerr}), see Fig.~\ref{fig:s2n}.
\end{itemize}
Fig.~\ref{fig:s2n} shows also the astrometric error, in arcsec, for
angular velocity $300$ arcsec/s (broken line).

\section{Correlation and orbit determination}
\label{sec:orbdet}

Given two or more observations sets, the main problem is to
identify which data belong to the same object (correlation
problem). Thus the orbit determination problem needs to be solved in
two stages: first different sets of observations need to be
correlated, then an orbit can be determined.

\subsection{Observations and attributables}

When a LEO is observed optically, even a comparatively short exposure
results in a trail. The detection of the trail, with the method
discussed in Sec.~\ref{sec:s2n}, is followed by the astrometric
reduction of the two ends of the trail. When the trail is too short to
show measurable curvature \cite{zoran05}, the information contained in
these two data points can be summarized in an attributable:
\[ 
{\cal A}=(\alpha,\delta,\dot\alpha,\dot\delta) \in [0,2\pi)
\times (\pi/2,\pi/2) \times \mathbb R^2\, ,
\]
representing the angular position and velocity of the body at a time
$t$. Usually $\alpha$ is the right ascension and $\delta$ the
declination with respect to an equatorial reference system (e.g.,
J2000).

To compute a 6 parameters orbit, we need 2 further
quantities. The values of range $\rho$ and range rate
$\dot\rho$ are not measured.  These two quantities, together with
${\cal A}$, allow us to compute the the Cartesian position and
velocity $(\mathbf r,\dot{\mathbf r})$ at time $\bar t=t-\rho/c$:
\begin{equation}\label{posvel}
\mathbf r=\mathbf q+\rho\hat{\bm\rho}\ \ ,\ \ \dot{\mathbf
  r}=\dot{\mathbf
  q}+\dot\rho\hat{\bm\rho}+\rho\frac{d\hat{\bm\rho}}{dt}\ \ \
, \ \ \
\frac{d\hat{\bm\rho}}{dt}=\dot\alpha\hat{\bm\rho}_\alpha+
\dot\delta\hat{\bm\rho}_\delta\ ,
\end{equation}
\[
\hat{\bm\rho}=(\cos\alpha\cos\delta,\sin\alpha\cos\delta,\sin\delta)\,
,\ \ \
\hat{\bm\rho}_\alpha=\partial\hat{\bm\rho}/\partial\alpha\,
,\ \ \
\hat{\bm\rho}_\delta=\partial\hat{\bm\rho}/\partial\alpha\, .
\]
In the above formulae the observer position $\mathbf q$ and velocity
$\dot{\mathbf q}$ are known functions of time.

\subsection{Keplerian integrals method}

In our simulations we used the Keplerian integrals method to produce
preliminary orbits from two attributables ${\cal A}_1$, ${\cal A}_2$
of the same object at two epoch times, $t_1$ and $t_2$. This is used
as the first, the most difficult, step of the correlation procedure.
The algorithm for the asteroid case is introduced in \cite{gronchi},
while in \cite{farnocchia} it is adapted to space debris.  The method
has already been successfully applied to the orbit determination of
GEOs in \cite{GEO}.  We recall here the basic steps.

Using (\ref{posvel}), the 2-body energy ${\cal E}$ and the angular
momentum $\mathbf c$ can be expressed as function of ${\cal A}$,
$\rho$ and $\dot\rho$. If the orbit between $t_1$ and $t_2$
is well approximated by a Keplerian one, we have:
\begin{equation}\label{keplerian_system}
\begin{cases}
{\cal E}({\cal A}_1,\rho_1,\dot\rho_1)-
{\cal E}({\cal A}_2,\rho_2,\dot\rho_2)=0\\
\mathbf c({\cal A}_1,\rho_1,\dot\rho_1)-
\mathbf c({\cal A}_2,\rho_2,\dot\rho_2)=0
\end{cases}
\ \ \ \rho_1>0\ ,\ \rho_2>0\, .
\end{equation}
From a solution of the above system we obtain 2 sets of
orbital elements at times $\bar t_1$ and $\bar t_2$:
\[
(a_j,e_j,I_j,\Omega_j,\omega_j,\ell_j)\ ,\ j=1,2\ .
\]
where $a$ is the semimajor axis, $e$ the eccentricity, $\Omega$ the
longitude of node, $\omega$ the argument of perigee, and $\ell$ the
mean anomaly. The first four Keplerian elements
$(a_j,e_j,I_j,\Omega_j)$ are functions of the 2-body energy and
angular momentum vectors ${\cal E}_j$, $\mathbf c_j$, and are the same
for $j=1,2$. In addition, there are compatibility conditions to be
satisfied if the two attributables belong to the same object:
\begin{equation}\label{compatibility_conditions}
\omega_1=\omega_2\ ,\ \ell_1=\ell_2+n(\bar t_1-\bar t_2)\ ,
\end{equation}
where $n=n(a)$ is the mean motion.
The algorithm also provides a covariance matrix for the computed orbit
and a $\chi$ value to measure the discrepancy in the compatibility
conditions. If such $\chi$ is smaller than some maximum value
$\chi_{max}$ the orbit is accepted.

\subsection{Precession model}
To deal with LEOs it is necessary to
generalize the method, including the effect due to the Earth
oblateness. The secular equations for Delaunay's variables $\ell$,
$\omega$, $\Omega$, $L=\sqrt{\mu a}$, $G=L\sqrt{1-e^2}$ and $Z=G\cos
I$ are (e.g., \cite{roy}[Sec. 11.4]):
\begin{equation}\label{delaunay_eq}
\begin{cases}
  \displaystyle\bar{\dot\ell}=n-\frac{3}{4}n
  \left(\frac{R_\oplus}{a}\right)^2\frac{J_2(1-3\cos^2I)}
  {(1-e^2)^{3/2}}\\
  \stackrel{}{\displaystyle\bar{\dot
      \omega}=\frac{3}{4}n\left(\frac{R_\oplus}
      {a}\right)^2\frac{J_2(4-5\sin^2I)}{(1-e^2)^2}}\\
  \stackrel{}{\displaystyle\bar{\dot
      \Omega}=-\frac{3}{2}n\left(\frac{R_\oplus}
      {a}\right)^2\frac{J_2\cos I}{(1-e^2)^2}}\\
  \stackrel{}{\bar{\dot L}=\bar{\dot G}=\bar{\dot Z}=0}
\end{cases}\ ,
\end{equation}
where $J_2$ is the coefficient of the second zonal spherical harmonic
of the Earth gravity field and $R_\oplus$ is the Earth radius. In this
case we cannot use the conservation of angular momentum, since it
precedes up to $5^\circ\cos I$ /day. Thus we consider the parametric
problem obtained by imposing $\dot\Omega=K$, with $K$ constant:
\begin{equation}\label{parametric_system}
\begin{cases}
  {\cal E}_1-{\cal E}_2=0\\
  R\,\mathbf c_1-R^T\,\mathbf c_2=0
\end{cases}
\end{equation}
where $R$ is the rotation by $\Delta\Omega/2=K(\bar t_2-\bar t_1)/2$
around $\hat{\mathbf z}$. Thus for a fixed value of $K$ the problem
has the same algebraic structure of the unperturbed one.  The
compatibility conditions contain the perigee precession and the
secular perturbation in mean anomaly, related to the one of the node
by
\[
\omega_1=\omega_2+K C_\omega\,(\bar t_1-\bar t_2)\, ,\ \
\ell_1=\ell_2+(n+K C_\ell)(\bar t_1-\bar t_2)\, ,
\]
where $C_\omega, C_\ell$ can be deduced from
(\ref{delaunay_eq}). Thus we can compute the $\chi(K)$.

We set up a fixed point iterative procedure, defined as follows:
\begin{enumerate}
\item given the orbital elements $E_{i}$ at step $i$ we compute the
  corresponding value of $K$ by (\ref{delaunay_eq});
\item for the computed value of $K$ we solve the corresponding system
  (\ref{parametric_system});
\item we select the solution with lowest
  $\chi$ and compute the corresponding orbital elements $E_{i+1}$;
\item we start again from 1, until convergence.
\end{enumerate}
To find a starting guess we select a value $K_0$ either using the
circular orbits corresponding to each of the two attributables (see
\cite{fujimoto}) or by assuming $K_0=0$. Among the obtained solutions,
we cannot know which one could lead to convergence. Since
the value of $K_0$ is quite arbitrary, it is not safe to select only
those with a low value of $\chi$. Thus we use all the solutions as
possible starting guess for the iterative procedure.

\subsection{Correlation confirmation}
\label{sec:correl_confirm}

The orbits obtained using the methods described in
the previous sections are just preliminary orbits, solution of a
2-body approximation, or possibly of a $J_2$-only problem. They have
to be replaced by least squares orbits, with a dynamical model
including all the relevant perturbations.

Even after confirmation by least squares fit, some correlations can be
\emph{false}, that is the two attributables might belong to different
objects. This is confirmed by the tests with real data reported in
\cite{GEO} and \cite{gronchi}. Thus every linkage of two attributables
needs to be confirmed by correlating further attributables.
The process of looking for a third attributable which can also be
correlated to the other two is called attribution
\cite{milani_recovery,milani_attrib}. From the available
2-attributable orbit with covariance we predict the attributable
${\cal A}_P$ at the time $t_3$ of a third attributable ${\cal A}_3$,
and test the compatibility with ${\cal A}_3$. For the successful
attributions we proceed to the differential corrections.

The procedure is recursive, that is we can use the 3-attributable
orbit to search for attribution of a fourth attributable, and so
on. This generates a very large number of many-attributable orbits,
but there are many duplications, corresponding to adding them in a
different order.  Thus a normalization step is needed to remove
duplicates (e.g., $A=B=C$ and $A=C=B$) and inferior correlations
(e.g., $A=B=C$ is superior to both $A=B$ and to $C=D$, thus both are
removed). This procedure, called \textit{correlation management}, is
important to reduce the risk of false correlation, because they are
eliminated by superior and discordant true correlations
\cite{milani_id}.  As an alternative, we may try to merge two
discordant correlations (with some attributables in common), by
computing a common orbit.  For a description of the sequence of steps
in this procedure see \cite{zoran05}.

\subsection{Updating well constrained orbits}
\label{sub:numbering}

When an orbit is very well constrained and the number of involved
trails is high enough, the possibility for the corresponding
correlation to be false can be ruled out. In this case the
corresponding object is cataloged and \textit{numbered} (that is, the
object is given a unique identifier). There are no fixed rules for
space debris to establish whether an object can be numbered. Our (safe
but arbitrary) criterion was the following: if the correlation
involves more than 10 trails the object is numbered. This criterion
makes sense, since 10 trails, with an observing strategy giving
only one exposure per pass, corresponds to many revolutions of the
object.

When an object is numbered, or anyway when the number of trails is
significant, it is possible to use the orbit with its covariance as a
basis for attribution of other trails to the same object. This
procedure should be done before attempting to find previously unknown
objects in the data, to reduce the data set to which the most
computationally complex procedure is applied. This is especially
necessary because the brightest objects are observed much more often
than the others. If the procedure of attribution of trails to known
objects is successful, after few weeks from the beginning of the
catalog build-up the leftover trails to be correlated are
comparatively few.

Note that the procedure outlined in this section essentially amounts
to a simulation not just of the performance of the sensors network,
but also of the operations of an orbit determination data center. The
performance of this element of a space awareness system is by no means
less important than the performance of the observing hardware.

\section{Results of the Catalog build-up simulation}
\label{sec:survey_res}

The first simulation was the one of the catalog build-up. Namely, we
assumed to start with no information on the satellite/debris
population in the region of interest (high LEO) and attempted to build
up an orbit catalog, to be compared with the MASTER one. For this
purpose we conducted a large scale simulation, with two
different population samples, dubbed Population 1 and 2.

\subsection{Population samples}
\label{sub:samples}

Population 1 was made up of 912 objects, randomly selected among those
with $8 < d <27$ cm, and $1300\leq h_p<2000$ km. Population 2 included
1104 objects, with $5<d<25$ cm, and $1000 \le h_p < 1300$ km. Objects
larger than 27 cm and 25 cm respectively were not included in the
simulation, since the amount of observational data allows us to attain
a complete catalog in a comparatively short time span. This is
confirmed not only by some preliminary tests but also by the analysis
of this simulation, which shows that all the objects larger than 20 cm
are cataloged within 1 month.

These choices were driven by the hypothetical performances of an
enhanced radar sensor, to be used along with the optical ones in the
ESA-SSA program. Figure~\ref{fig:simul12_2months} shows two rectangles
on the perigee altitude-diameter plane, which enclose the regions
corresponding to the selected population samples.  For the considered
values of the perigee altitude, we have drawn two curves, representing
the observing capabilities of a \textit{baseline radar} and of an
\textit{enhanced radar}. The corresponding assumptions on the radar
systems are rather arbitrary, thus these performance curves have to be
considered just as a comparison benchmark: indeed, the purpose of the
investigation was to determine the requirement for a future radar
system, and to see how much this could change by assuming the
cooperation of a network of optical sensors. These curves give the
minimum diameter $d_{min}$ of observable objects as function of
$h_{p}$, according to the following law:
\begin{equation}
d_{ min}=\frac{h_{p}^2}{h_{ref}^2}d_{ref}\,,
\label{eq:radarcurves}
\end{equation}
where $h_{ref}$ and $d_{ref}$ are reference values for the
perigee altitude and the diameter. The equation (\ref{eq:radarcurves})
is a simple consequence of the inverse fourth power dependence of
radar $S/N$ from distance. We set $h_{ref}=2000$ km and $d_{\rm
  ref}=32$ cm for the baseline radar. For the enhanced radar we took
$h_{ref}=2000$ km and $d_{ref}=20$ cm.
The two rectangles containing the populations used in the simulation
cover a region including the performance curve of the enhanced radar.

To analyze the results, we split the population in three
orbital classes:
\begin{itemize}
\item[] {\bf LEO}: LEO resident objects, with semimajor axis $a \le
  R_{\oplus} + $ 2000 km;
\item[] {\bf PLEO}: ``partial LEO'', with perigee inside but apogee
  beyond the LEO region: $h_p< 2000$ km, $a > R_{\oplus} + 2000$ km,
  and limited eccentricity, so that $a < 25000$ km;
\item[] {\bf HLEO}: LEO transit objects, with $h_p < 2000$ km and
  large eccentricity, so that $a \ge 25000$ km. This class includes in
  particular GTO and Molniya objects.
\end{itemize}
The HLEO objects spend only a small fraction of their orbital period
below $2000$ km of altitude, and have a large velocity when at
perigee, thus an angular velocity with respect to the observing
station larger than the resident LEO. Thus they are much more
difficult to be observed by the network selected for resident
LEO. They could be the target of a survey with a completely different
observing strategy. E.g., the GTO orbits can be observed with longer
exposures and non-sidereal tracking while near apogee \cite{GTO}. Thus
the results in our simulations about HLEO are not very significant:
these objects should be considered in the context of a simulation with
an observing strategy adapted to higher orbits, such as MEO.

In the MASTER population model, smaller objects typically are sampled,
i.e. a single object represents many fragments. Simple de-sampling
(e.g., by assigning different mean anomalies to each fragment, keeping
all the other orbital elements the same) could affect our simulations
by generating false correlations. Therefore the sampled fragments were
treated as single objects in the simulation process. A \emph{sampling
  factor} was associated to each object to represent the number of
fragments within the same orbit in the final analysis of the results.

\subsection{Results}
\label{sub:survey_res}

\begin{figure}[h!]
\centering
\includegraphics[width=0.9\textwidth]{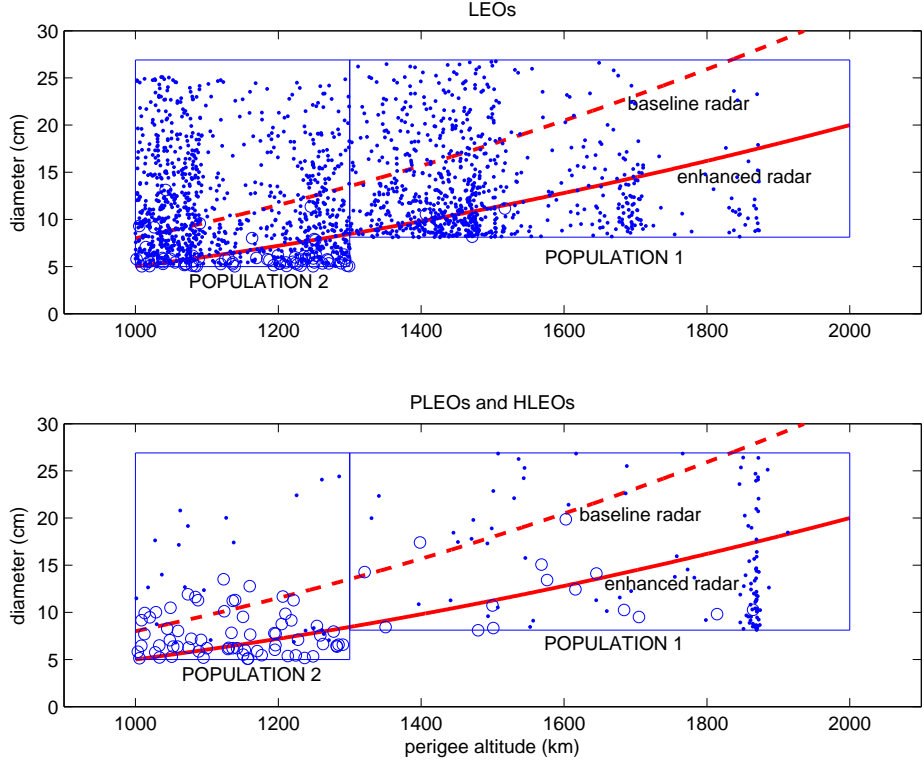}
\caption{Results of catalog build-up after 2 months}
\label{fig:simul12_2months}
\end{figure}

Figure~\ref{fig:simul12_2months} shows the results of both population
samples in the perigee altitude-diameter plane after 2 months of
survey observations. The rectangles enclose the regions of the plane
corresponding to the selected population samples. The parabolic curves
represent two possible radar performances in the case of a baseline
radar (dashed line) and of an enhanced one (solid line). The top
figure shows the results for resident LEOs, while the bottom one is
for transit objects (PLEOs and HLEOs). Points indicate successful
orbit determination ($\geq$ 3 trails involved), while circles
correspond to failure. The reasons for the occurrence of this last
possibility are essentially three:
\begin{enumerate}
\item lack of observations;
\item the observations were very few and too distant in time, so that
  correlation was not possible;
\item the observations had low accuracy, so that even when the
  correlation of two trails was successful, the orbit accuracy was not
  enough for the attribution of further trails.
\end{enumerate}
For LEOs, there are many failures only in the low part of the region
corresponding to Population 2. This means that this possibility can
happen only for some very small objects, with diameter $< 8$ cm.

Tables~\ref{tab:survey_res1} and \ref{tab:survey_res2} detail the
results for each population sample. The factor associated to the
sampled objects is not taken into account. This means that a group of
fragments with the same orbit is treated as a single object. In the
tables are reported in this order the total number of objects in the
sample, the number of cataloged objects, the number of objects
observed but not cataloged, and the number of never observed
objects. The row in parenthesis contains the number of objects, among
those with some observations but without orbit, for which a total of
only 1-2 trails were available in 2 months of survey.

\begin{table}
  \caption{\label{tab:survey_res1} Results of catalog build-up after 2 
    months, Population 1}
\centering
\begin{tabular}{|l|c|c|c|c|}
\hline
{\bf POP 1}&Total & LEO & PLEO & HLEO\\
\hline
No. Objects &912 &796 &97 &19\\
Orbits Computed &894& 793& 95& 6\\
Obj. without orbit &10&3&2&5\\
(with 1-2 Tr.) &(3)&(0)&(0)&(3)\\
Obj. not observed &8 &0 &0 &8\\
\hline
\end{tabular}
\end{table}

\begin{table}
\caption{\label{tab:survey_res2} Results of catalog build-up after 2 months, Population 2}
\centering
\begin{tabular}{|l|c|c|c|c|}
\hline
{\bf POP 2}&Total & LEO & PLEO & HLEO\\
\hline
No. Objects &1104& 1014& 62& 28\\
Orbits Computed &965 &942 &21& 2\\
Obj. without orbit&92&67&16&9\\
(with 1-2 Tr.)&(13)&(9)&(2)&(2)\\
Obj. not observed &47& 5& 25& 17\\
\hline
\end{tabular}
\end{table}

To analyze the results, we measured the efficiency of the catalog
build-up, defined as the ratio between the number of reliable orbits
computed (with at least 3 trails) and the total number of sampled
objects. We measured also the efficiency only for those objects lying
in the region of the perigee altitude-diameter plane above the curve
of the enhanced radar. These efficiencies were computed accounting for
the number of clones for the sampled objects (sampling factor), to
count both the successes and the failures of correlation. The detailed
results are shown in Tables~\ref{tab:eff1} and \ref{tab:eff2}.

\begin{table}
  \caption{\label{tab:eff1} Efficiency of catalog build-up after 2 
    months, Population 1}
\centering
\begin{tabular}{|l|c|c|c|c|}
\hline
{\bf POP 1}&Total &LEO &PLEO &HLEO\\
\hline
Eff. Catalog &98.1\%& 99.6\%& 97.9\%& 31.6\%\\
Eff. above Radar &98.6\%& 99.8\%& 97.2\%& 71.4\%\\
\hline
\end{tabular}
\end{table}

\begin{table}
  \caption{\label{tab:eff2} Efficiency of catalog build-up after 2 months, Population 2}
\centering
\begin{tabular}{|c|c|c|c|c|}
\hline
{\bf POP 2}&Total &LEO &PLEO &HLEO\\
\hline
Eff. Catalog &82.8\%& 86.6\% &37.5\%& 7.1\%\\
Eff. above Radar& 93.7\% &98.9\% &97.2\% &25.0\%\\
\hline
\end{tabular}
\end{table}

The results for Population 1 were overwhelming. For $h_p >
1300$ km, we computed the orbits of almost all the objects above 8 cm
($98.1\%$) and of $99.6\%$ of resident LEOs. If we consider only the
objects above the curve of enhanced radar performance in
Fig.~\ref{fig:simul12_2months}, then the total efficiency is $98.6\%$,
while for resident LEOs it is $99.8\%$. A few objects in highly
eccentric transit orbits (e.g., GTO) were lost, but they should be
observed with a different strategy, as outlined previously. 

For Population 2 the total efficiency is $82.8\%$,
while considering only resident LEOs it is $86.6\%$; if we take into
account only the region above the curve of the enhanced radar the
corresponding percentages are $93.7\%$ and $98.9\%$. Moreover, the
outcome of the two simulations suggests that the problem is not the
altitude but the diameter: indeed, even for Population 2, among the
objects larger than 8 cm, $95.8\%$ had a reliable orbit, and for
resident LEOs the percentage was $99.0\%$. Smaller debris require a
larger telescope, whatever the height.

To get the big picture of the results, we should take into account
also the biggest objects of the MASTER population. Moreover, the
results of Population 1 have to be weighted twice, since we randomly
selected only half of the objects. By assuming that all the biggest
objects were cataloged after 2 months of survey, we get an overall
efficiency 97.9\% for the objects above the enhanced radar curve. By
considering only LEOs the efficiency is
99.0\%. Figure~\ref{fig:simul12_2months} shows that most of the
difficulties in the catalog build-up concern the region of orbital
perigee altitude less than $1100$ km. If we discard this region, the
total efficiency reaches 98.9\%, and for LEOs 99.8\%.

\section{Orbit improvement phase}
\label{sec:task_results}

Tasking observations are possible if, by moving the telescope at
non-sidereal rate, the object remains in the telescope field of
view. This is achieved if the angular error in prediction is small
enough. 
Moreover, it is desirable to see the object with most signal
concentrated on a single pixel, thus obtaining a higher signal with
respect to the survey mode.

\subsection{Feasibility of tasking observations}
\label{sec:task_assumptions}

We considered the numbered objects obtained after 1 month of catalog
build-up and we propagated the orbits for 1 week after the last
observation. Then the angular error in prediction and the relative
error in angular velocity were computed, by means of the covariance
matrices. 

Denoting by $\Gamma_{P}$ and $\Gamma_{V}$ the marginal covariance
matrices for position and velocity respectively, and by
$\lambda_{P}, \lambda_{V}$ their maximum eigenvalues, we used
the following upper bounds for the angular error and the relative
error in velocity:
\[
\Delta(\alpha,\delta)\le \frac{\sqrt{\lambda_{P}}}{h_{p}}\,,
\quad
\frac{|\Delta {\bf v}|}{|{\bf v}|}\le 
\frac{\sqrt{\lambda_{V}}}{|{\bf v}|}\ .
\]
These approximations are justified for LEOs by the fact that they move
on almost circular orbits, so that the radial component of the
velocity is small. On the other hand high eccentricity transit objects
are observed near the perigee.

Taking into account that the angular velocity of the objects
in our sample was less than 2000 arcsec/s, to see the object in a
single pixel during tasking, the relative error in angular velocity
had to be less than $7.5 \times 10^{-4}$. This was achieved for all
the numbered objects. Moreover the angular error in position was
always less than 95 arcsec. Then the possibility of follow up was
confirmed and we were allowed to take $T=1$ (no trailing loss) in the
$S/N$ formula, with a gain in $S/N$ for faint objects by a factor at least
14.

Besides assuming $T=1$, we supposed the scheduler to be capable of
selecting 1/6 of the observation period, corresponding to the best
phase angle, thus gaining more than 1 magnitude with respect to survey
observations. The combination of these two effects increases the
number of observations, because a trail with a very low $S/N$ in
survey mode becomes detectable in tasking mode as a point source; in
addition, in many cases the astrometric accuracy is significantly
improved.

To perform the follow up simulation, we chose the same two populations
used for the survey step, but we excluded HLEOs, since the observation
strategy was not suitable for this kind of objects, as confirmed by the
results of the catalog build-up.
At the starting time of the simulation, the sampled objects were
assumed to be cataloged with sufficiently well determined orbits.
Thus, we assumed that the switch to tasking mode occurred for all the
cataloged orbits at the same time. In reality, the decision to include
an object in the tasking mode scheduling should be performed on an
individual basis.
The same procedure of the survey simulations was adopted, meaning that
the previously computed catalog was not used and the orbit
determination was performed again from scratch. The difference was in
the kind and the amount of data, as follow up observations were more
numerous and accurate than the survey ones.

\subsection{Accuracy envelope norm}

The orbit improvement goal was defined by choosing an accuracy
envelope to which compare the accuracy of the computed orbits. First
we fixed the object-centered reference frame
$\{\mathbf{u},\mathbf{v},\mathbf{w}\}$, where $\mathbf{u}$ is the
direction from the Earth center to the object center, $\mathbf{w}$ is the
direction of the angular momentum vector and $\mathbf{v} = \mathbf{w}
\times \mathbf{u}$.  Then the accuracy ellipsoid in position for
resident LEOs was defined by the quantities $4, 30, 20$ m along the directions
$\mathbf{u},\mathbf{v},\mathbf{w}$ respectively, for PLEOs by
 $10, 60, 200$ m. The accuracy ellipsoid in velocity was $20, 4, 20$
mm/s both for resident LEOs and PLEOs. These ellipsoids were
derived from collision avoidance requirements. 

To see if the improved orbits were inside the envelope, we had to verify
that the confidence ellipsoids, representing the uncertainty of the
orbits in position and velocity, were fully contained in the ellipsoids
defined by the accuracy envelope. 

Let $\CC_{req}$ be the marginal normal matrix associated to the
accuracy envelope for position (for LEOs, $\CC_{req}=diag [1/4^2,
1/30^2, 1/20^2]$ in the OCRF), and let $\CC_{conf}$ be the marginal
normal matrix of the orbit. Then the two ellipsoids relative to
position are defined respectively by the equations
\[ 
\xx^T \CC_{req} \xx=1\,, \quad \xx^T \CC_{conf}
\xx=1\,.
\] 
To satisfy the envelope, the maximum of the function $\xx^T \CC_{req} \xx$,
as $\xx$ varies in the confidence ellipsoid of the orbit, must be less
than 1.
By Lagrange multiplier theorem a necessary and sufficient condition
for this is that $\lambda \le 1$, for any $\lambda$ for which
$\det(\CC_{req}-\lambda \CC_{conf})=0$.

Let $\Cov_{req}=\CC_{req}^{-1}$ and $\Cov_{conf}=\CC_{\rm
  conf}^{-1}$ denote the marginal covariance matrices and let
$\mathbf{P}_{req}$ a matrix such that $\CC_{\rm
  req}=\mathbf{P}_{req}^T\mathbf{P}_{req}$ (it exists,
because $\CC_{req}$ is symmetric and non-negative). Then the values
of $\lambda$ such that $\det(\CC_{req}-\lambda \CC_{conf})=0$
are the eigenvalues of the matrix $\mathbf{P}_{req}\Cov_{\rm
  conf}\mathbf{P}_{req}^T$.
Thus we defined the {\bf accuracy envelope norm} in position as
the square root of the maximum eigenvalue of the
matrix $\mathbf{P}_{req}\Cov_{conf}\mathbf{P}_{req}^T$.
For the velocity we used a completely analogous definition.

\subsection{Results} 

The simulation covered three weeks of tasking. All the numbered
orbits were propagated for 1 week after the last observation and
compared to the accuracy envelope. The comparison was performed
through the accuracy envelope norm previously defined, which
is $\le 1$ whenever the uncertainty is smaller than the envelope. 

Tables~\ref{tab:improvleo1} and \ref{tab:improvleo3} summarize the
results for LEOs after the first and the third week respectively.  The
values reported are: percentage of objects with accuracy envelope norm
$\le 1$ in both position and velocity, maximum norms achieved, maximum
angular error in position, percentage of objects with angular error
$\le 1.5$ arcsec (the pixel size). 
For PLEOs the results are not statistically representative, because
in both simulations they were less then 100 objects.

The conclusion of the orbit improvement simulation is that just one
week of data is not enough, even with the higher observation frequency
and accuracy of the tasking phase, while 3 weeks are more than enough.

\begin{table}
\caption{\label{tab:improvleo1} Results for LEOs after 1 week of tasking}
\centering
\begin{tabular}{|l|c|c|c|c|c|}
\hline
{\bf LEO}& norms & max norm &
max norm &max ang.& ang. err. \\
1 WEEK&$\le 1$&position&velocity&err. (arcsec)&$\le 1.5$ arcsec\\
\hline
POP. 1&76.8\%&5.69&5.19&12.96&79.2\%\\
POP. 2&68.6\%&6.18&6.5&24.12&63.1\%\\
\hline
\end{tabular}
\end{table}

\begin{table}
\caption{\label{tab:improvleo3} Accuracy of the results for LEOs 
  after 3 weeks of tasking}
\centering
\begin{tabular}{|l|c|c|c|c|c|}
\hline
{\bf LEO}& norms & max norm &
max norm &max ang.& ang. err. \\
3 WEEKS&$\le 1$&position&velocity&err. (arcsec)&$\le 1.5$ arcsec\\
\hline
POP. 1&99.5\%&1.25&1.09&2.88&98.2\%\\
POP. 2&99.9\%& 1.03&0.97& 2.88&94.3\%\\
\hline
\end{tabular}
\end{table}

\begin{figure}[h!]
\centering
\includegraphics[width=0.9\textwidth]{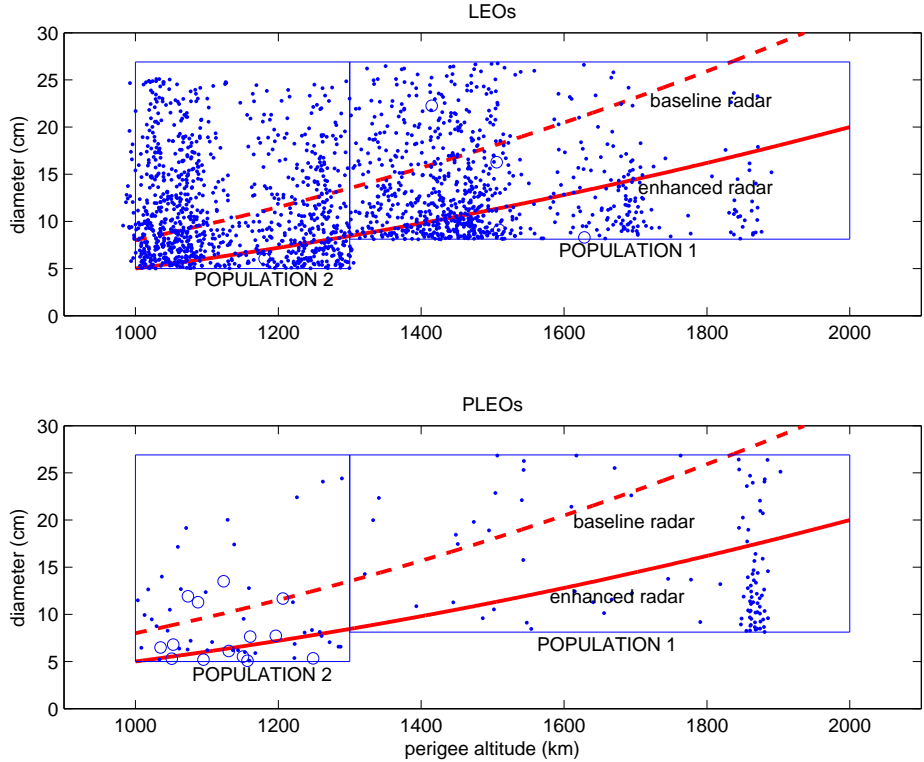}
\caption{Results of both simulations after 3
  weeks of tasking}
\label{fig:improv3w}
\end{figure}

Figure~\ref{fig:improv3w} shows the results after 3 weeks. The top
figure is for resident LEOs while the bottom one is for PLEOs. Points
indicate that both the accuracy norms in position and velocity are
$\le 1$, while circles mean that one of the above norms is $> 1$. It turns
out that we succeeded for most of LEOs, while there were a few problems
for low altitude PLEOs.

Finally, it must be stressed that the relative error in velocity
prediction was always less than $10^{-4}$, which means that with
tasking observations there is no trailing loss.

\section{Fragmentation detection}
\label{sec:frag}

A particularly demanding situation for a surveillance network occurs
whenever a fragmentation happens in Earth orbit. The resulting cloud
of fragments can pose significant immediate or short term risks for
other space assets in the same orbital region
\cite{rossi1999}. Therefore two ad hoc simulations were performed to
verify the capability of our proposed optical network to detect and
characterize a fragmentation event in high LEO.  In general terms, the
detection of fragmentation implies to detect a significant number of
fragments. Another possibility is to detect the change in the orbital
elements of the target object due to the energy of the event. Note
that, in the case of an active satellite, if no other information is
available, it could be difficult to discriminate between a change due
to a maneuver and a collision.

We simulated the fragmentation of a satellite due to both an explosion
and a catastrophic collision. To produce the cloud of fragments we
used, for both cases, the ISTI implementation of the NASA-JSC model
\cite{nasa_break}.

Given a mass $M_T$ of the target object, a cloud of fragments was
generated with the proper size distribution given by the NASA
model. According to the model, each fragment was characterized by the
following elements: size $d$, mass, area, $\Delta v$.  From the whole
cloud of fragments, we selected those objects with $d >10$ cm and
$\Delta v<100$ m/s, that is we concentrated on the core of the cloud.
Given the Cartesian coordinates of the parent body, the isotropic
$\Delta v$ relative to each fragment was added to these coordinates,
thus obtaining the state vector of each fragment.  From the Cartesian
coordinates, the orbital elements of each fragment were then computed.
The orbital elements of the fragments were then propagated for 21 days
and the simulated observations were computed. The simulations
assumed observations in survey mode, with trailing loss. 

First we simulated the explosion of a satellite with mass $M_T = 1000$
kg. The satellite was supposed to be on a circular orbit at 1400 km of
altitude, with an inclination of $74^\circ$. The model produced a
cloud of 230 fragments larger than 10 cm. Out of them, 175 fragments
had $\Delta v< 100$ m/s (and 8 fragments had $\Delta v < 10$
m/s). The detection of the fragments by the whole network of 7 ground
stations was simulated.  During the first day after the explosion it
was not possible to detect any debris, while, in the second and in the
third day, about half of the debris were detected (26.3 \% in the
second day, 57.7 \% in the third day and already 98.9 \% in the fourth
day). This situation depends on the fact that in the first day the
fragments cloud had still a non-homogeneous distribution around the
Earth. The debris created by the explosion formed a cloud concentrated
around the explosion point. Therefore, in the first hours after the
event, the detection of the concentrated cloud was similar to the
detection of a single object. If the cloud was not passing in the right
moment of the night above a station with favorable meteorological conditions,
its detection was not possible. On the contrary, a few days after
the explosion, all the debris were detected because they had a
homogeneous distribution in a torus around the parent body orbit. It
is therefore important to remind that the debris detection after a
fragmentation depends on the ground station meteorological data.  In
particular, in the first day only 2 ground stations were able to
observe debris, but the debris cloud was not passing over these
stations.

In the second simulation, the collision of a debris with mass 10 kg,
against a satellite with mass $M_T = 1000$ kg, with an impact velocity
of $v = 9$ km/s, was simulated.  In the NASA model, a collision is
considered catastrophic whenever the ratio between the kinetic energy
of the projectile and the mass of the target exceeds $40\,000$
J/kg. This condition is largely satisfied by our event, so we are
considering a catastrophic collision. The target satellite was
supposed to be on the same orbit as in the explosion case.  The model
produced a cloud of 998 fragments larger than 10 cm. Out of them, 393
fragments had $\Delta v< 100$ m/s (and 7 fragments had $\Delta v <
10$ m/s).  The evolution of the cloud of fragments was similar to the
one discussed in the explosion case.

The detection statistics was similar to the one for the explosion case
(33.6 \% of detected objects in the second day, 57.7 in the third day,
66.4 in the fourth day and 99.7 \% in the fifth day) and the same
considerations about the meteorological conditions apply here. In the
first day after the collision only 2 ground stations were able to
observe some debris.

After simulating the observations performed by the ground network, we
started the correlation and orbit determination process, to check how
soon after the event it is possible to have robust information on the
fragmentation.  Clearly, if 2 days were enough to observe a quite high
percentage of fragments, they were not enough to compute reliable
orbits for the observed fragments. There were too few observations and
in fact no orbits were available after only two days for both the
simulations. After 4 days the situation improved and 90 orbits for the
collision fragments and 46 orbits for the explosion fragments were
computed, if we considered as acceptable correlations of at least 3
trails and discarded the ones between only 2 trails. By comparison with
the ground truth (i.e., with the original fragments generated in the
simulations) we found that some of the correlations were false, that is
they put together observations of different objects. Anyway, all the
false correlations had only 3 trails and this was true for both the
simulations in the entire period examined.  Note that, in the case of
the fragments clouds, it is quite natural to have false correlations
even among 3 trails, because all the fragments have very similar
orbits. To exclude false correlations, in this peculiar case, we
considered only orbits fitting $\geq$ 4 trails.

In our analysis we considered as reliable the orbits which fit at
least 5 trails and as numbered the orbits with at least 10 trails. We
found that, after only 2 weeks, all the objects were numbered both in
the explosion and in the collision case. Table~\ref{tab:frag} gives the
summary of the orbit determination process, by showing the number of
orbits computed with at least 5 and 10 trails, for both clouds in the
two weeks following the fragmentation events.

\begin{figure}[h]
\centering
\includegraphics[width=0.7\textwidth]{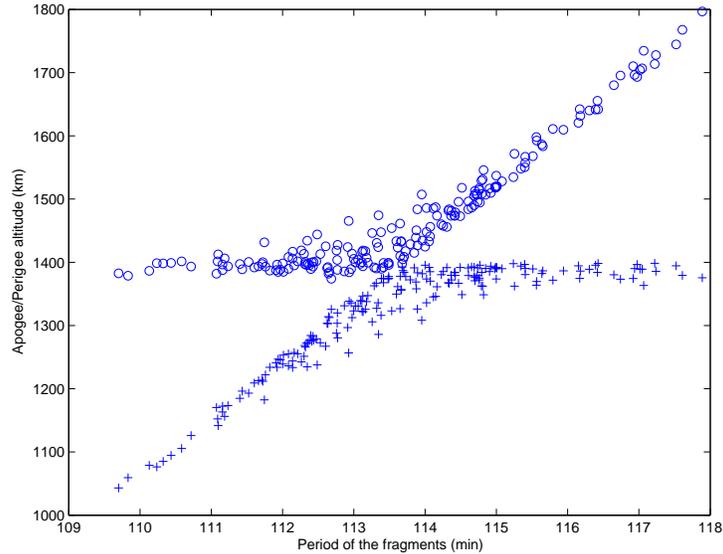}
\caption{Gabbard diagram of the collision fragments cataloged
within 6 days after the event}
\label{fig:gabbard}
\end{figure}

The graphical tool commonly used to characterize an in-orbit
fragmentation is the Gabbard diagram, plotting the apogee/perigee
height of the fragments as a function of their orbital
period. Figure~\ref{fig:gabbard} shows the Gabbard diagram (circles
for apogee, crosses for perigee) for the collision case, 6 days after
the event. The X-shape, a typical signature of a fragmentation, is
already visible.

\begin{table}
\caption{\label{tab:frag} Computation of the orbits of fragments for 
  both explosion and collision}
\centering
\begin{tabular}{|c|c|c|c|c|}
\hline
& \multicolumn{2}{c|}{Explosion}  &\multicolumn{2}{c|}{Collision} \\
\hline
Days after&Orbits & Orbits &Orbits & Orbits \\
fragmentation event & $\ge 5$ trail & $\ge 10$ trails & $\ge 5$ trail 
& $\ge 10$ trails\\
\hline
2  & 0	   &   0   & 0    &    0 \\
4  & 26	   &   0   & 63   &    0 \\
6  & 99	   &   62  & 192  &  141 \\
8  & 139   &   124 & 307  &  287 \\
10 & 154   &   149 & 333  &  329 \\
12 & 169   &   169 & 366  &  366 \\
14 & 175   &   175 & 392  &  392 \\
16 & 175   &   175 & 393  &  393 \\
\hline
\end{tabular}
\end{table}

In conclusion, our study suggests that the simulated network
of telescopes can detect and catalog the fragments generated by a
catastrophic event in high LEO after just a few days from the event.
Some caveats and conclusions can be stated.  First, the detection of a
stream of fragments, with low ejection velocity, within 24 hours is
like detecting a single object, because the fragments are not spread
along the entire orbit. For this reason it might not be possible to
detect a large fraction of the fragments in 1 day, if bad meteorological
conditions are present on critical stations.  Then, the Gabbard
diagram, built with the output of the orbit determination simulation
after 6 days, shows that the orbital information is more than enough
to assess the fragmentation event (parent body, energy, etc.).


\section{Conclusions}
\label{sec:conclusion}
The results of the catalog build-up simulation show that more than
98\% of the LEO objects with perigee height above 1100 km and diameter
greater than 8 cm can be cataloged in 2 months. As
Fig.\ref{fig:simul12_2months} shows, a central area around 1100 km of
orbital perigee altitude has been identified where radar sensors and the
optical network should operate in a cooperative way. All the numbered
orbits are accurate enough to allow follow up observations with no
trailing loss, and the orbit accuracy from the improved orbits is
compliant with the accuracy we have used, which corresponds to
collision avoidance requirements. Finally the simulated network of
telescopes is able to detect and catalog the fragments generated by a
catastrophic event just a few days after the event.

The significance of our results for the design of a Space Situational
Awareness system is in the possibility to use optical sensors to
catalog and follow up space debris on orbits significantly lower than
those previously considered suitable. Of course this is true only
provided a list of technological assumptions, both in hardware and in
software, spelled out in Section~\ref{sec:assumptions}, are
satisfied. If these technologies are available, then it is possible to
trade off between an upgraded system of optical sensors and a radar
system with higher energy density.

\section*{Aknowledgments}
We wish to thank CGS for providing us with the expected sensor
performances and realistic simulated observations, taking into account
the characteristics of the Fly-Eye telescope, the meteorological model
to account for cloud cover, and the statistical model for the signal
to noise ratio. This work was performed under ESA/ESOC contract
22750/09/D/HK \textit{SARA-Part I Feasibility study of an innovative system
for debris surveillance in LEO regime}, with CGS as prime contractor.





\bibliographystyle{elsarticle-num}



\end{document}